\input ./style/arxiv-general.cfg
\documentclass[MSNbibl,number,citesort,secthm,seceqn,dvips]{arxbj}
\makeatletter
   \@ifpackageloaded{graphicx}{}{\usepackage{graphicx}}
\makeatother

\volume{22}
\issue{3}
\pubyear{2016}
\firstpage{1729}
\lastpage{1744}
\doi{10.3150/15-BEJ709}
\docsubty{FLA}

\makeatletter
\newtheorem{lemma}{Lemma}[section]
\newremark{rem}{Remark}[section]
\newremark{exa}{Example}[section]
\makeatother

\begin{document}
\begin{frontmatter}

\title{Optimal classification and nonparametric regression for
functional data}
\runtitle{Optimal classification and nonparametric regression for
functional data}

\begin{aug}
\author[A]{\inits{A.}\fnms{Alexander}~\snm{Meister}\corref{}\ead[label=e1]{alexander.meister@uni-rostock.de}}
\address[A]{Institut f\"ur Mathematik, Universit\"at Rostock, D-18051
Rostock, Germany.\\ \printead{e1}}
\end{aug}

%
\received{\smonth{4} \syear{2014}}
%
\revised{\smonth{12} \syear{2014}}

\begin{abstract}
We establish minimax convergence rates for classification of functional
data and for nonparametric regression with functional design variables.
The optimal rates are of logarithmic type under smoothness constraints
on the functional density and the regression mapping, respectively.
These asymptotic properties are attainable by conventional kernel
procedures. The bandwidth selector does not require knowledge of the
smoothness level of the target mapping. In this work, the functional
data are considered as realisations of random variables which take
their values in a general Polish metric space. We impose certain metric
entropy constraints on this space; but no algebraic properties are required.
\end{abstract}

\end{frontmatter}

\section{Introduction} \label{1}

In many statistical applications, the empirical data cannot be
described by random vectors in a Euclidean space $\mathbb{R}^d$. Still
one can often reasonably define a distance between the possible
realisations of the observations. Then parts of the data are supposed
to take their values in a non-empty Polish metric space $({\mathcal X},\rho
)$ where the corresponding probability measure has the corresponding
Borel $\sigma$-field $\mathfrak{B}({\mathcal X})$ as its domain. Note that
a separable and complete metric space is called a Polish metric space.

Within that general framework, the analysis of functional data has
attained increasing attention (see, e.g., the book of   Ramsay and Silverman \cite{r22} for an introduction to the topic). Therein ${\mathcal
X}$ denotes some appropriate function space, for example, the set of
all continuous and bounded functions on $[0,1]$ or the set of all
measurable and squared-integrable functions on that domain. The current
work is mainly motivated by this research field; whereas, in general,
the elements of ${\mathcal X}$ are not imposed to be functions or
equivalence classes of functions, which opens up new perspectives for
extensions to even more complex types of data. In particular, we use
only topological properties of the set ${\mathcal X}$; but no algebraic
structure on ${\mathcal X}$ is required (e.g., linear space, group, ring,
etc.). Therefore, tools from principal component analysis
(e.g., Benko,  H{\"a}rdle and Kneip \cite{r2}) or manifold representation (e.g., Chen and M{\"u}ller \cite{r9}) cannot be applied in this setting. Instead, we use arguments
based on covering and packing numbers from the approximation theory.
Such techniques are frequently used in empirical process theory in
order to study the parameter set of a statistical experiment, which
consists of functions in nonparametric statistics (e.g., van~der Vaart and Wellner \cite{r25}, van~de Geer \cite{r24}). Contrarily, they have been
applied to learn about the sample set in quite few papers.

We focus on two widely studied problems in functional data analysis:
nonparametric regression (Section~\ref{2}) and classification (Section~\ref{3}). A review on existing literature is provided in the
corresponding sections. While a huge amount of literature is available
on these topics, only little has been known about the aspect of
asymptotic optimality of statistical procedures when the sample size
$n$ tends to infinity. The current note intends to advance the
understanding of those problems by providing the minimax convergence
rates for the statistical risks. The proofs are deferred to
Section~\ref{4}. Section~\ref{1.5} provides some essential topological tools which
are used in both Sections~\ref{2} and~\ref{3}.

\section{Entropy condition} \label{1.5}

In the following, we recall two concepts from approximation theory
(e.g., van~der Vaart and Wellner \cite{r25}, page 83, Definition~2.1.5 and
page 98, Definition~2.2.3): by ${\mathcal N}_{\mathcal X}(\delta,{\mathcal Y},\rho)$
we denote the covering number of some set ${\mathcal Y} \subseteq{\mathcal
X}$, that is, the minimal number of open $\rho$-balls in ${\mathcal X}$ with
the radius $\delta$ so that ${\mathcal Y}$ is a subset of the union of
these balls. If we stipulate in addition that the centers of those
balls lie in ${\mathcal Y}$, we call this quantity the intrinsic covering
number ${\mathcal N}_{\mathcal Y}(\delta,{\mathcal Y},\rho)$. The packing number
${\mathcal D}(\delta,{\mathcal Y},\rho)$ of the set ${\mathcal Y}$ describes the
maximal cardinality of a subset of ${\mathcal Y}$ such that $\rho
(x,y)>\delta$ for all elements $x\neq y$ of this subset. Also, we learn
from Kolmogorov and Tihomirov \cite{r19} and van~der Vaart and Wellner \cite{r25}, page~98, that
\begin{equation}
\label{eq2.1}
\mathcal{N}_{\mathcal Y}(\delta,{\mathcal Y},\rho) \leq {
\mathcal D}(\delta ,{\mathcal Y},\rho) \leq {\mathcal N}_{\mathcal Y}(\delta/2,{
\mathcal Y},\rho)\qquad \forall\delta>0.
\end{equation}
Also, we easily derive that
\begin{equation}
\label{eq2.2}
{\mathcal N}_{\mathcal X}(\delta,{\mathcal Y},\rho) \leq {
\mathcal N}_{\mathcal
Y}(\delta,{\mathcal Y},\rho) \leq {\mathcal
N}_{\mathcal X}(\delta/2,{\mathcal Y},\rho) \qquad \forall\delta>0.
\end{equation}
Now we classify a type of sets ${\mathcal Y}$ by their metric
entropy, which we define by
\[
\Phi(s,{\mathcal Y},\rho) := \log{\mathcal N}_{\mathcal X}(s,{\mathcal Y},
\rho) \qquad \forall s>0.
\]
Concretely, we assume that
\begin{equation}
\label{eq2.3} c_{x,0} s^{-\gamma} \leq \Phi(s,{\mathcal Y},\rho)
\leq c_{x,1} s^{-\gamma}
\qquad \forall s\in(0,s_0),
\end{equation}
for some fixed constants $s_0>0$, $0 < c_{x,0} < c_{x,1}$ and
$\gamma>0$. We write $B_{\mathcal Y}(x,r)   :=   \{y\in{\mathcal Y}
\dvt
\rho(x,y) < r\}$ for $x\in{\mathcal X}$ and $r>0$. We easily see that
$B_{\mathcal Y}(x,r) \in\mathfrak{B}({\mathcal X})$ for all $x\in{\mathcal X}$,
$r>0$ and ${\mathcal Y} \in\mathfrak{B}({\mathcal X})$. Condition (\ref
{eq2.3}) can be justified in many applications. Let us consider two
examples of classes ${\mathcal Y}$ which satisfy this condition.

\begin{exa}[(Classes of smooth functions)]\label{R.1}
We assume that our functional data $X_1,\ldots,X_n$ are located in a
class of smooth functions almost surely. We write $\lceil\alpha\rceil
$ for the smallest integer which is larger or equal to $\alpha>0$.
Precisely, we impose the H\"older constraints that ${\mathcal Y}$ consists
of functions $f$ mapping from $[0,1]^d$ to $\mathbb{R}$ such that all
partial derivatives of $f$ up to the order $\lceil\alpha\rceil- 1$
are bounded by a constant $M$; and that the $(\lceil\alpha\rceil-
1)$th partial derivatives satisfy the H\"older condition with the
exponent $\alpha-\lceil\alpha\rceil+ 1$ and again the constant $M$.
Also, we put ${\mathcal X} = C_0([0,1]^d)$ and $\rho$ equal to the supremum metric.

We learn from Theorem~2.7.1, page~155 in van~der Vaart and Wellner \cite{r25} that the upper bound in condition (\ref{eq2.3}) is satisfied
with $\gamma= d/\alpha$. Also, the corresponding lower bound can be
verified (see Kolmogorov and Tihomirov \cite{r19}).

Moreover, for any $\alpha>0$, the H\"older class ${\mathcal Y}$ is
relatively compact with respect to the supremum metric thanks to the
Arzel{\`a}--Ascoli theorem, from what follows compactness of the
closure~$\overline{{\mathcal Y}}$. We deduce that
\[
\lim_{\delta' \downarrow\delta} {\mathcal N}_{\mathcal X} \bigl(
\delta',\overline {{\mathcal Y}},\rho \bigr) \leq {\mathcal
N}_{\mathcal X}(\delta,{\mathcal Y},\rho) \leq {\mathcal N}_{\mathcal X}(
\delta,\overline{{\mathcal Y}},\rho) \qquad \forall\delta>0,
\]
since, for any cover of ${\mathcal Y}$ by the union of finitely
many open balls, the union of the corresponding closed balls covers
$\overline{{\mathcal Y}}$ (and so does the union of the corresponding\vspace*{1pt} open
balls with arbitrarily enlarged radius). Therefore, condition (\ref
{eq2.3}) is extended from ${\mathcal Y}$ to $\overline{{\mathcal Y}}$; and the
role of ${\mathcal Y}$ can be taken over by its closure.
\end{exa}

In general, the technique of the last paragraph in Example~\ref{R.1},
that is, switching to the closure of ${\mathcal Y}$, can be used to impose
without loss of generality that ${\mathcal Y}$ is closed -- and hence,
${\mathcal Y} \in\mathfrak{B}({\mathcal X})$~-- without any loss of generality
when condition (\ref{eq2.3}) is assumed.

\begin{exa}[(Classes of monotonic functions)] \label{R.2}
Now we consider the example of componentwise monotonic mappings from
the cube $[0,1]^d$ to $[0,1]$. The collection of these functions is
denoted by ${\mathcal Y}$. As the corresponding Polish metric space, we
choose ${\mathcal X} = L_p([0,1]^d)$, $p\geq1$, that is, the Banach space
of all Borel measurable functions $f$ from $[0,1]^d$ to $\mathbb{R}$
which satisfies $\int|f(x)|^p \,\mathrm{d}x < \infty$. Clearly, $\rho$ is the
metric generated by the $L_p([0,1]^d)$-norm.

Then Theorem~1.1 in Gao and Wellner \cite{r17} yields that ${\mathcal Y}$
satisfies condition (\ref{eq2.3}) with $\gamma= \max\{d,(d-1)p\}$ for
$d\geq2$ and $(d-1)p\neq d$. In the univariate setting $d=1$, the
upper bound part of condition (\ref{eq2.3}) with $\gamma=1$ follows
from Theorem~2.7.5 in van~der Vaart and Wellner \cite{r25}, page 159.
Therein we use that the covering number is bounded from above by the
bracketing number with doubled radius for the $L_p([0,1]^d)$-metric
$\rho$ (see page~84, van~der Vaart and Wellner \cite{r25}). On the other
hand, the according lower bound can be established by Proposition~2.1
in Gao and Wellner \cite{r17}.
\end{exa}

The following lemma provides a useful result for the upper bound
proofs in the following two sections.

\begin{lemma} \label{Ltop}
Let $({\mathcal X},\rho)$ be a Polish metric space. Take some ${\mathcal Y}\in
\mathfrak{B}({\mathcal X})$ which satisfies (\ref{eq2.3}), and let $P$ be
any probability measure on $\mathfrak{B}({\mathcal X})$ with $P({\mathcal
Y})=1$. We set
\[
\psi(x,h) := P \bigl(B_{\mathcal Y}(x,h) \bigr),\qquad h>0.
\]
Then we have
\[
P \bigl( \bigl\{x \in{\mathcal Y} \dvt \psi(x,h) \leq\delta \bigr\} \bigr) \leq
\delta\exp \bigl(c_{x,1} 4^\gamma h^{-\gamma} \bigr),
\]
for all $\delta>0$.
\end{lemma}

\section{Nonparametric regression} \label{2}

We observe the data set ${\mathcal Z}_n = \{(X_1,Y_1),\ldots,(X_n,Y_n)\}$
where the $X_j$ are i.i.d. random variables taking their values in the
Polish metric space $({\mathcal X},\rho)$ equipped with the corresponding
Borel $\sigma$-algebra $\mathfrak{B}({\mathcal X})$. The $Y_j$ are defined by
\begin{equation}
\label{eqmodel}
Y_j = g(X_j) + \varepsilon_j,
\end{equation}
where $g$ denotes some Borel measurable mapping from ${\mathcal
X}$ to $\mathbb{R}$; and the $\varepsilon_j$ are real-valued random
variables which satisfy
\begin{equation}
\label{eq3.2}
E(\varepsilon_1\vert X_1) = 0,
\qquad \operatorname{var}(\varepsilon _1\vert X_1) \leq c_v,
\qquad \mbox{$P_X$-a.s.},
\end{equation}
for some uniform constant $c_v$ where $P_X$ denotes the
probability measure on $\mathfrak{B}({\mathcal X})$ which is generated by
$X_1$. The random variables $(X_1,\varepsilon_1),\ldots,(X_n,\varepsilon
_n)$ are assumed to be i.i.d. Moreover, we assume that $X_1\in{\mathcal
Y}$ holds almost surely for some subset ${\mathcal Y}\in\mathfrak{B}({\mathcal
X})$. Our goal is to estimate the regression function $g$ based on the
data set ${\mathcal Z}_n$.

As a usual condition in nonparametric regression, we impose some
smoothness constraints on the regression function $g$. Precisely, we
introduce the class ${\mathcal G} = {\mathcal G}_{\beta,C}$ of all Borel
measurable mappings $g$ from ${\mathcal X}$ to $\mathbb{R}$ such that $\sup_{y\in{\mathcal Y}} |g(y)| \leq C$ and
\[
\bigl|g(y)-g(z) \bigr| \leq C \rho(y,z)^\beta \qquad \forall y,z\in{\mathcal Y}
\]
with $C>0$, $\beta\in(0,1]$. Critically, we remark that our
framework is restricted to smoothness degrees $\beta$ which are smaller
or equal to one. An extension to higher smoothness levels seems
difficult as ${\mathcal X}$ is not equipped with any algebraic structure so
that no common definitions of Taylor series can be applied. Approaches
to local linear methods, which should capture all smoothness levels
smaller than two, are provided in
Berlinet,
Elamine and Mas \cite{r3} and Mas \cite{r20}; while, in these papers, ${\mathcal X}$ is assumed to be a Hilbert
space -- transferred to our notation.

Whereas linear models for $g$ (along with generalizations) are popular
in functional regression problems (e.g., Hall and Horowitz \cite{r18},
Meister \cite{r21}), fully nonparametric approaches to the regression
function have also received considerable attention. We refer to the
book of Ferraty and Vieu \cite{r14} for a comprehensive review on kernel
methods for functional covariates. In Ferraty \textit{et al.} \cite{r15}, a generic
upper bound is derived for the uniform rate of convergence. Recently,
Forzani, Fraiman and Llop
\cite{r16} consider consistency of nonparametric functional
regression estimation in the setting of a metric space without any
imposed algebraic structure. In a similar setting, Biau,
C{\'e}rou and Guyader \cite{r5}
establish upper bounds on an integrated risk for the convergence rates
of the functional $k$-nearest neighbor estimator when $\beta= 1$ (in
our notation). The convergence rates used in that paper are of
logarithmic type. However, minimax optimality is apparently not studied
in this work.

To our best knowledge, the only approach to rate-optimal nonparametric
functional regression estimation is given by Mas \cite{r20}, who uses
principal component analysis on ${\mathcal X}$ and specific conditions on
these components. The attained rates are faster than any logarithmic
rates but slower than any polynomial rate in the non-Gaussian case. In
our setting where the design distribution obeys the condition (\ref
{eq2.3}), the minimax convergence rates are different.
We consider estimators $\hat{g}$ of $g$ which are Borel measurable
mappings from ${\mathcal X}^{n+1}$ to $\mathbb{R}$ and which are squared
integrable with respect to the design measure $P_X$ after inserting the
data, regardless of their realization. Also, we impose that $g \in
L_2(P_X)$, that is, the Hilbert space of all squared\vspace*{1pt} integrable and
measurable functions with respect to $P_X$. Then we are guaranteed that
$\|\hat{g}(\cdot,{\mathcal Z}_n) - g\|_{P_X}^2$ is a real-valued random
variable where $\|\cdot\|_{P_X}$ denotes the $L_2(P_X)$-norm.

We take the Nadaraya--Watson estimator for functional data,
\begin{equation}
\label{eqFV}
\hat{g}(x) :=
\cases{\displaystyle \hat{A}(x) /
\hat{B}(x), & \quad$\mbox{if }\hat{B}(x) > \delta_n$, \vspace*{3pt}
\cr
0, & \quad\mbox{otherwise},}
\end{equation}
where
\begin{eqnarray*}
\hat{A}(x) & := & \frac{1}n \sum_{j=1}^n
Y_j K_h \bigl(\rho(x,X_j) \bigr),
\\
\hat{B}(x) & := & \frac{1}n \sum_{j=1}^n
K_h \bigl(\rho(x,X_j) \bigr).
\end{eqnarray*}
However, we have modified the concept by adding the
truncation to the denominator $\hat{B}(x)$ where the ridge parameter
$\delta_n>0$ remains to be selected. Moreover, $h>0$ denotes a
bandwidth parameter and $K \dvtx \mathbb{R}\to\mathbb{R}$ a kernel function.
We employ the notation $K_h := K(\cdot/h)$ (without dividing by $h$).
For simplicity, we choose that $K = 1_{[0,1)}$. We provide the
following asymptotic result.
\begin{thm} \label{T2}
Let ${\mathcal Y}\in\mathfrak{B}({\mathcal X})$ such that \textup{(\ref{eq2.3})} holds
true. We consider model \textup{(\ref{eqmodel})} under the condition \textup{(\ref{eq3.2})}.
Then, for any sequence $\{P_{X,n}\}_n$ of design measures on
$\mathfrak{B}({\mathcal X})$ with $P_{X,n}({\mathcal Y}) = 1$ for all $n$, the
estimator $\hat{g}$ in (\ref{eqFV}) satisfies
\[
\sup_{g\in{\mathcal G}} \int E
\bigl|\hat{g}(x) - g(x)\bigr|^2
\,\mathrm{d}P_X(x) = {\mathcal O} \bigl(\{\log n\}^{-2\beta/\gamma} \bigr),
\]
under\vspace*{1pt} the kernel choice $K = 1_{[0,1)}$ and the parameter
selection $\delta_n = n^{-\eta}$, $\eta\in(0,1/2)$ and $h = \{d \log
n\}^{-1/\gamma}$ with $d\in (0,\eta c_{x,1}^{-1} 4^{-\gamma} )$.
\end{thm}

\begin{rem}
Under the additional assumption $P_X \in{\mathcal R}_X$, which says that
\[
P_X \bigl(B_{\mathcal Y}(y,\delta) \bigr) \geq c_{x,3}
\delta\exp \bigl(-c_{x,4} \delta^{-\gamma} \bigr) \qquad \forall\delta
\in(0,1), y\in{\mathcal Y},
\]
which can be shown to be non-empty for some positive constants
$c_{x,3}$ and $c_{x,4}$ and any compact~${\mathcal Y}$, we can also derive
the following upper bound on the pointwise risk:
\[
\sup_{P_X \in{\mathcal R}_X} \sup_{g\in{\mathcal G}} \sup
_{x\in{\mathcal Y}} E\bigl|\hat{g}(x) - g(x)\bigr|^2 = {\mathcal O} \bigl(
\{\log n\}^{-2\beta/\gamma
} \bigr),
\]
under the same conditions on $K$, $\delta_n$ and $h$ as in Theorem~\ref
{T2} except that $d\in(0,\eta/c_{x,4})$.
\end{rem}

We consider model (\ref{eqmodel}) with the additional condition that
the $\varepsilon_j$ are i.i.d. random variables with a continuously
differentiable density function $f_\varepsilon$ with finite Fisher
information, that is,
\begin{equation}
\label{eqFisher}
\int\bigl|f_\varepsilon'(x)\bigr|^2 /
f_\varepsilon(x)\, \mathrm{d}x < \infty.
\end{equation}
Moreover, all the $X_1,\varepsilon_1,\ldots,X_n,\varepsilon
_n$ are independent. Also we impose compactness of the set ${\mathcal Y}$
from Theorem~\ref{T1}. The following theorem provides an asymptotic
lower bound for the estimation of $g$ with respect to the pointwise
estimation error as well as the integrated risk.
\begin{thm} \label{T1}
Let ${\mathcal Y}\in\mathfrak{B}({\mathcal X})$ be compact and assume that
\textup{(\ref{eq2.3})} holds true. We consider model~\textup{(\ref{eqmodel})} under
independent additive regression errors $\varepsilon_j$, $j=1,\ldots,n$
with a density $f_\varepsilon$ which satisfies \textup{(\ref{eqFisher})}.
\begin{longlist}[(a)]
\item[(a)] Then there exists a sequence of design measures $P_{X,n}$ on
$\mathfrak{B}({\mathcal X})$ with $P_{X,n}({\mathcal Y}) = 1$ for all~$n$, such
that no sequence of estimators $\{\hat{g}_n\}_n$ based on the data
${\mathcal Z}_n$ satisfies
\[
\sup_{g\in{\mathcal G}} \int E\bigl|\hat{g}(x) - g(x)\bigr|^2
\,\mathrm{d}P_{X,n}(x) = \mathrm{o} \bigl(\{\log n\}^{-2\beta/\gamma} \bigr).
\]
\item[(b)] For any sequence of design measures $P_{X,n}$ on
$\mathfrak{B}({\mathcal X})$ with $P_{X,n}({\mathcal Y}) = 1$ for all $n$ and
for any sequence of estimators $\{\hat{g}_n\}_n$ based on the data
${\mathcal Z}_n$, we have
\[
\liminf_{n\to\infty} \sup_{g\in{\mathcal G}} \sup
_{x\in{\mathcal Y}} P \bigl[\bigl|\hat{g}(x) - g(x)\bigr|^2 > c\cdot\{\log
n\}^{-2\beta/\gamma} \bigr] > 0,
\]
for some constant $c$ depending on $C$ and $\beta$.
\end{longlist}
\end{thm}

Theorem~\ref{T1} establishes minimax optimality of the convergence
rate attained in Theorem~\ref{T2} in two views. Part (a) shows that
there exists a sequence of design measures such that the integrated
risk does not converge with faster rates. Obviously, we cannot obtain
such a result for any design measure: if $P_{X,n}$ was a one-point
measure then just the average of the $Y_1,\ldots,Y_n$ would be a
consistent estimator with the usual parametric rate. In part (b), we
prove that no matter what the design measure looks like, one is not
able to obtain faster pointwise convergence rates simultaneously for
all $x\in \mathcal{Y}$, even with respect to the weak rates.

An important and widely studied issue in nonparametric regression is
bandwidth selection. The minimax convergence rates are of slow
logarithmic type. However, the bandwidth selector in Theorem~\ref{T2}
leads to the optimal rates while it can be used without knowing the
smoothness degree $\beta$. This selector is fully deterministic, which
means that no data-driven procedure (e.g., cross validation, Lepski's
method, etc.) is required in order to achieve the optimal convergence
rates. It is remarkable that the same effects occur in nonparametric
deconvolution from supersmooth error distributions (see, e.g., Fan \cite{r13}) and other severely ill-posed inverse problems. We face a
bias-dominating problem, that is, the variance term is asymptotically
negligible under the optimal bandwidth selection. In other
bias-dominating problems, sharp asymptotics have been studied (Butucea
and Tsybakov \cite{r7}). It is an interesting question for future research
if those results apply to the current problem as well.

\section{Classification} \label{3}

The problem of classifying functional data has also stimulated great
research activity (e.g., Ferraty and Vieu \cite{r14},
Carroll, Delaigle and Hall \cite{r8},
Delaigle and Hall \cite{r10,r11}, Biau, Bunea and
Wegkamp  \cite{r4}). It has its
applications in the fields of biometrics, genetics, recognition of
sounds, technometrics, etc. Classification problems are closely linked
to the field of statistical learning theory (e.g., Vapnik \cite{r26}). We
choose the model of supervised classification. Concretely, we observe
some random variable $Z$ taking its values in some Polish metric space
${\mathcal X}$ -- and in ${\mathcal Y} \in\mathfrak{B}({\mathcal X})$ almost
surely. We assume that we have two groups $0$ and $1$ and our goal is
to decide whether $Z$ should be categorized as a member of group $0$ or
$1$. The groups $0$ and $1$ are characterized by the probability
measures $P_X$ and $P_Y$ on $\mathfrak{B}({\mathcal X})$, respectively. One
does not know these measures; however, i.i.d. a training sample
$(Z_j,W_j)$, $j=1,\ldots,n$ is available where the $W_j$ are binary
random variables and $W_j=b$, $b=0,1$, indicates that $Z_j$ has the
probability measure $P_X$ and $P_Y$, respectively. Moreover, $Z$ is
independent of all training data.

In order to specify all admitted probability measures $P_X$ and $P_Y$,
we impose that
\begin{eqnarray}
 (P_X,P_Y) \in {\mathcal P}_{\kappa}
&:=  & \bigl\{(P,Q) \dvt \mbox{$P$ and $Q$ are probability measures on $
\mathfrak{B}({ \mathcal X})$ so that}
\nonumber
\\[-8pt]
\label{eqcl1}
\\[-8pt]
\nonumber
&& \hspace*{19pt}\qquad P({\mathcal Y}) = Q({\mathcal Y}) = 1, \operatorname{TV}(P,Q) \geq
\kappa \bigr\},
\end{eqnarray}
for some $\kappa> 0$ where $\operatorname{TV}(P,Q)$ denotes the
total variation distance between some measures $P$ and $Q$,
\[
\operatorname{TV}(P,Q) := \sup_{A\in\mathfrak{B}({\mathcal X})} \bigl|P(A) - Q(A) \bigr|.
\]
With respect to the set ${\mathcal Y,}$ we assume condition (\ref
{eq2.3}).

Unlike in classification problems for data in $\mathbb{R}^d$, $d\in
\mathbb{N}$, we face the problem that no spatially homogeneous measure
(e.g., Lebesgue--Borel measure, Haar measure) exists on $\mathfrak
{B}({\mathcal X})$ so that no density of $P_X$ and $P_Y$ can be defined
with respect to such a measure. Nevertheless, $P_X$ and $P_Y$ are
dominated by their sum measure $Q := P_X+P_Y$. We write $p_X$ and $p_Y$
for the Radon--Nikodym derivatives $p_X := \mathrm{d}P_X/\mathrm{d}Q$ and $p_Y := \mathrm{d}P_Y/\mathrm{d}Q
= 1-p_X$. We impose some smoothness constraints on both $p_X$ and $p_Y$ via
\begin{eqnarray}
(P_X,P_Y) \in {\mathcal P}_{C,\beta,\kappa}
&:= &  \bigl\{(P_X,P_Y) \in {\mathcal P}_\kappa
\dvt \exists{\mathcal Y}_0\in\mathcal{P}({\mathcal Y})\cap\mathfrak
{B}({ \mathcal X}) \mbox{ with }[P_X+P_Y]({\mathcal
Y}_0)=2 \mbox{ s.t.}
\nonumber
\\[-8pt]
\label{eqcl2}
\\[-8pt]
\nonumber
 &&\hspace*{50pt}\qquad \bigl|p_X(y) - p_X(z) \bigr| \leq C
\rho^\beta(y,z), \forall y,z \in{\mathcal Y}_0 \bigr\},
\end{eqnarray}
with $C>0$ and $\beta\in(0,1]$ -- analogously as in Section~\ref{2} in the regression setting. Therein $\mathcal{P}({\mathcal Y})$
denotes the power set of ${\mathcal Y}$.

A (supervised) classifier $\varphi$ is defined as a Borel measurable
mapping from ${\mathcal X}^{n}\times\{0,1\}^n \times{\mathcal X}$ to $\{0,1\}
$. Clearly, the sample $(Z_1,\ldots,Z_n,W_1,\ldots,W_n,Z)$ is inserted
into $\varphi$ and $\varphi=b$, $b=0,1$, means categorizing $Z$ as a
member of group $b$. We define the excess risk of classification by
\[
{\mathcal E}_{n}(\varphi) := \sup_{(P_X,P_Y) \in{\mathcal P}_{C,\beta
,\kappa}}
\bigl(P_{X,Y,X} [\varphi= 1 ] + P_{X,Y,Y} [\varphi= 0 ] - 1 +
\operatorname{TV}(P_X,P_Y) \bigr),
\]
in order to evaluate the accuracy of some classifier $\varphi$.
The excess risk is the sum of the probabilities of misclassification
into group $0$ and $1$, respectively, reduced by $1 - \operatorname{TV}(P_X,P_Y)$. Therein $P_{X,Y,X}$ and $P_{X,Y,Y}$ indicate that $Z$
has the probability measure $P_X$ or $P_Y$, respectively. It is well
known that the excess risk of the Bayes classifier
\[
\varphi_B({\mathbf z},{\mathbf w},z) := %
\cases{ 0, &
\quad$\mbox{if }p_X(z) \geq1/2$, \vspace*{3pt}
\cr
1, &\quad
\mbox{otherwise},}
\]
vanishes if ${\mathcal P}_{C,\beta,\kappa}$ was replaced by some
two-element set $\{P_X,P_Y\}$, that is, if $P_X$ and $P_Y$ were known.

Our goal is to find a classifier $\varphi$ which minimizes the excess
risk asymptotically as $n,m$ tend to infinity. To our best knowledge
optimal convergence rates for classification of functional data have
been unexplored so far; whereas for finite-dimensional data they have
been studied, for example, in Yang \cite{r27,r28} and Audibert and Tsybakov \cite{r1}. Considering the Bayes classifier, it is reasonable to mimic the
unknown densities $p_X$ and $p_Y$ by some appropriate estimators based
on the data $Z_1,W_1,\ldots,Z_n,W_n$ (also see, e.g.,
Biau, Bunea and Wegkamp \cite{r4} or Ferraty and Vieu \cite{r14}). We employ the classifier
\begin{equation}
\label{eqclest}
\varphi(Z_1,\ldots,Z_n,W_1,
\ldots,W_n,Z) = \cases{ 0, & \quad$\mbox{if } \hat{p}_X(Z)
\geq\hat{p}_Y(Z)$,\vspace*{3pt}
\cr
1, & \quad \mbox{otherwise,}}
\end{equation}
where
\begin{eqnarray*}
\hat{p}_X(z) & := & \sum_{j=1}^n
(1-W_j) \cdot K \bigl(\rho (z,Z_j)/h \bigr) \Big/ \sum
_{j=1}^n (1-W_j),
\\
\hat{p}_Y(z) & :=& \sum_{j=1}^n
W_j \cdot K \bigl(\rho(z,Z_j)/h \bigr) \Big/ \sum
_{j=1}^n W_j,
\end{eqnarray*}
if $\sum_{j=1}^n (1-W_j) \in(0,n)$; otherwise put $\hat{p}_X(z)=0$ or
$\hat{p}_Y(z)=0$ by convention. Therein we apply some kernel $K$ and
bandwidth parameter $h>0$ as in Section~\ref{2}. We stipulate that
enough data $Z_j$ from both $P_X$ and $P_Y$ are available; concretely,
we impose
\begin{equation}
\label{eqnm}
P[W_1=1]=w \qquad \mbox{for some fixed value $w\in(0,1)$}.
\end{equation}
The asymptotic performance of the classifier (\ref{eqclest})
is studied in the following theorem.
\begin{thm} \label{T3}
We consider the model of supervised classification. Let ${\mathcal Y}\in
\mathfrak{B}({\mathcal X})$ such that \textup{(\ref{eq2.3})} holds true. Moreover,
we assume \textup{(\ref{eqnm})}. Then the excess risk of the classifier $\varphi
$ in \textup{(\ref{eqclest})} attains the following uniform upper bound:
\[
{\mathcal E}_{n}(\varphi) = {\mathcal O} \bigl((\log
n)^{-\beta/\gamma} \bigr),
\]
under the kernel choice and the bandwidth selection from
Theorem~\ref{T2}.
\end{thm}

While Theorem~\ref{T3} can be proved directly, it follows from Theorem~\ref{T2} by the general argument that the excess mass is bounded from
above by the integrated squared regression risk (see, e.g.,
Devroye, Gy{\"o}rfi and Lugosi \cite{r12}, page~104). Furthermore, we mention that, in the setting of
Theorem~\ref{T3}, we could relax the assumptions contained in ${\mathcal
P}_{C,\beta,\kappa}$ to $\kappa=0$. Still the condition $\kappa>0$ is
realistic as $P_X$ and $P_Y$ should not become too close to each other;
otherwise, the classification problem makes no sense. Finally, in
Theorem~\ref{T4} we will establish optimality of the convergence rates
from Theorem~\ref{T3} with respect to an arbitrary sequence of classifiers.
\begin{thm} \label{T4}
We consider the model of supervised classification. Let ${\mathcal Y}\in
\mathfrak{B}({\mathcal X})$ such that \textup{(\ref{eq2.3})} holds true. Moreover,
we assume \textup{(\ref{eqnm})}. Fix some $\kappa>0$ sufficiently small (but
independent of~$n$). Let $\{\varphi_n\}_n$ be an arbitrary sequence of
(supervised) classifiers where $\varphi_n$ is based on the data
$(Z_1,\ldots,Z_{n},W_1,\ldots,W_{n},Z)$. Then we have
\[
\mathop{\liminf}_{N\to\infty} (\log n )^{\beta/\gamma} {\mathcal
E}_{n}(\varphi_n) > 0.
\]
\end{thm}

The optimal convergence rates in Theorems \ref{T3} and \ref{T4}
correspond to those established in Section~\ref{2} in the regression
problem. Note that there we consider the squared risk. Again, we
realize that the bandwidth selector in Theorem~\ref{T3} does not
require knowledge of the smoothness level $\beta$ and, still, it leads
to the optimal speed of convergence.

\section{Proofs}\vspace*{-12pt} \label{4}
\begin{pf*}{Proof of Lemma~\protect\ref{Ltop}}
Let $X$, $Y$ be some
independent random variables with the induced measure $P$. Note that
$\psi(X,h)$ can be viewed as the conditional expectation of the random
variable $1_{[0,h)}(\rho(X,Y))$ given $X$ so that the random mapping
$\psi(X,h)$ is measurable, thus a random variable. By the factorization
lemma of the conditional expectation, the mapping $x \mapsto\psi
(x,h)$, $x\in{\mathcal X}$, is measurable so that ${\mathcal Y}_{h,\delta} :=
 \{x\in{\mathcal Y} \dvt \psi(x,h) \leq\delta \}$ lies in $\mathfrak
{B}({\mathcal X})$. Furthermore, we obtain that
\begin{eqnarray*}
&& {\mathcal Y}_{h,\delta}  = \bigcup_{j=1}^{{\mathcal N}_{\mathcal
Y}(h/2,{\mathcal Y},\rho)}
\bigl\{y\in B_{\mathcal Y}(y_{j},h/2) \dvt \psi (y,h) \leq\delta
\bigr\},
\end{eqnarray*}
where $\{y_{1},\ldots,y_{{\mathcal N}_{\mathcal Y}(h/2,{\mathcal Y},\rho
)}\} \subseteq{\mathcal Y}$ denotes an intrinsic $h/2$-cover of ${\mathcal Y}$
with respect to the metric $\rho$. By $J$ we denote the collection of
all $j=1,\ldots,{\mathcal N}_{\mathcal Y}(h/2,{\mathcal Y},\rho)$ such that the set
$ \{y\in B_{\mathcal Y}(y_{j},h/2)   \dvt   \psi(y,h) \leq\delta \}$
is not empty. For any $j\in J$, there exists some $y\in B_{\mathcal
Y}(y_{j},h/2)$ with $P(B_{\mathcal Y}(y,h)) \leq\delta$. We have $B_{\mathcal
Y}(y_{j},h/2) \subseteq B_{\mathcal Y}(y,h)$ so that $P (B_{\mathcal
Y}(y_{j},h/2) ) \leq\delta$. We deduce that
\begin{eqnarray*}
&& P({\mathcal Y}_{h,\delta})  \leq \sum_{j\in J} P
\bigl(B_{\mathcal
Y}(y_{j},h/2) \bigr) \leq \delta{\mathcal
N}_{\mathcal Y}(h/2,{\mathcal Y},\rho ) \leq \delta\exp \bigl(c_{x,1}
4^\gamma h^{-\gamma} \bigr),
\end{eqnarray*}
when combining (\ref{eq2.2}) and (\ref{eq2.3}).
\end{pf*}

\begin{pf*}{Proof of Theorem~\protect\ref{T2}}
For any $g\in{\mathcal G}$ we
derive that
\begin{eqnarray}
\nonumber
&& E \bigl\{\bigl|\hat{g}(x) - g(x)\bigr|^2 \vert X_1,
\ldots,X_{n} \bigr\}
\\
\nonumber
&&\quad \leq 1 \bigl\{\hat{B}(x)>\delta_n \bigr\}
\hat{B}^{-2}(x) E \bigl\{ \bigl|\hat {A}(x) - g(x) \hat{B}(x)\bigr|^2
\vert X_1,\ldots,X_{n} \bigr\} + g^2(x) \cdot1
\bigl\{\hat{B}(x)\leq\delta_n \bigr\}\qquad\hspace*{7pt}
\nonumber
\\[-8pt]
\label{eqPr1}
\\[-8pt]
\nonumber
&&\quad \leq 2C^2 h^{2\beta} + 2 c_v \cdot1
\bigl\{\hat{B}(x)>\delta _n \bigr\} \hat{B}^{-2}(x)
n^{-2} \sum_{j=1}^n
1_{[0,h)} \bigl(\rho(X_j,x) \bigr) + C^2 \cdot1
\bigl\{\hat{B}(x)\leq\delta_n \bigr\}
\\
\nonumber
&&\quad \leq 2C^2 h^{2\beta} + 2 c_v
n^{-1} \delta_n^{-2} + C^2 \cdot1
\bigl\{\hat{B}(x)\leq\delta_n \bigr\},
\end{eqnarray}
holds almost surely under the convention $0\cdot\infty= 0$.
We realize that $E \hat{B}(x) = \psi_n(x,h) := P_{X,n}(B_{\mathcal
Y}(x,h))$. By the inequality,
\[
1 \bigl\{\hat{B}(x)\leq\delta_n \bigr\} \leq 1 \bigl\{\bigl|\hat{B}(x)-
\psi _n(y,h)\bigr| \geq\delta_n \bigr\} + 1 \bigl\{
\psi_n(x,h) \leq2\delta _n \bigr\},
\]
applying the expectation to both sides of (\ref{eqPr1})
leads to
\begin{eqnarray}
\label{eqPr2}
&&\hspace*{-12pt} E\bigl|\hat{g}(x) - g(x)\bigr|^2 \leq 2C^2
h^{2\beta} + 2 c_v n^{-1} \delta_n^{-2}
+ C^2 \cdot\delta_n^{-2} \operatorname{var} \hat {B}(x)
+ C^2 \cdot1 \bigl\{\psi_n(x,h) \leq2\delta_n
\bigr\},\qquad
\end{eqnarray}
where $\operatorname{var} \hat{B}(x) \leq n^{-1} \psi_n(x,h)$.
Putting $x = X_{n+1}$ (i.e., an independent copy of $X_1,\ldots,X_n$)
and applying the expectation to both sides of (\ref{eqPr2}) leads to
\[
E\bigl|\hat{g}(X_{n+1}) - g(X_{n+1})\bigr|^2 \leq
2C^2 h^{2\beta} + \bigl(2c_v + C^2
\bigr) n^{-1} \delta_n^{-2} + C^2 P
\bigl[\psi_n(X_{n+1},h) \leq2\delta _n \bigr].
\]
Putting $\psi_n = \psi$, $X_{n+1}=X$ and $2\delta_n = \delta
$, Lemma~\ref{Ltop} yields that
\[
P \bigl[\psi_n(X_{n+1},h) \leq2\delta_n \bigr]
\leq2 n^{4^\gamma c_{x,1}
d -\eta}.
\]
Due to the constraint on $d$ the term $2C^2 h^{2\beta}$ is
asymptotically dominating, which provides the desired upper bound on
the considered risk with uniform constants on $g\in{\mathcal G}$.
\end{pf*}

\begin{pf*}{Proof of Theorem~\protect\ref{T1}}
(a)  We introduce some
sequence $(\delta_n) \downarrow0$. As ${\mathcal Y}$ satisfies (\ref
{eq2.3}), the packing number has the lower bound
\[
{\mathcal D}(\delta_n,{\mathcal Y},\rho) \geq m_n :=
\exp \bigl(c_{x,0} \delta_n^{-\gamma} \bigr),
\]
due to (\ref{eq2.1}) and (\ref{eq2.2}). This implies the
existence of some $z_{1,n},\ldots,z_{m_n,n} \in{\mathcal Y}$ such that the
balls $B_{j,n}:= B_{\mathcal Y}(z_{j,n},\delta_n/4)$, $j=1,\ldots,m_n$ are
pairwise disjoint. This statement can be strengthened to the result
that the $\rho$-distance between the sets $B_{j,n}$ and $\bigcup_{k\neq
j} B_{k,n}$ is even bounded from below by $\delta_n/2$. We specify $P_X
= P_{X,n}$ as the discrete uniform distribution on the grid $\{
z_{1,n},\ldots,z_{m_n,n}\}$.

We use the function $\vartheta(t)   =   \exp\{1/(t^2 - 1)\} \cdot
1_{(-1,1)}(t)$, $t\in\mathbb{R}$. Thus, $\vartheta$ is differentiable
infinitely often on the whole real line, yielding that
\[
\bigl|\vartheta(t) - \vartheta(s)\bigr| \leq \min \bigl\{\|\vartheta\|_\infty ,
\bigl\|\vartheta'\bigr\|_\infty|t-s| \bigr\} \leq \max \bigl\{\|
\vartheta\| _\infty, \bigl\|\vartheta'\bigr\|_\infty \bigr\}
\cdot|t-s|^\beta.
\]
We construct the regression\vspace*{-3pt} curves
\[
g_\theta(x) = \sum_{j=1}^{m_n}
\theta_j d h_n^\beta\vartheta \bigl(\rho
(z_{j,n},x)/h_n \bigr),
\]
with the vector $\theta= (\theta_1,\ldots,\theta_{m_n}) \in
\{0,1\}^{m_n}$ and $h_n := \delta_n/4$, for some $d>0$. As $(h_n)_n$ is
bounded from above, the constraint $\sup_{g\in{\mathcal G}} \|g\|_\infty
\leq C$ can be satisfied by choosing $d>0$ small enough. For all
$y_1,y_2 \in{\mathcal Y}$ there exist at most one $j_1$ and one $j_2$ such
that $y_l \in B_{j_l,n}$, $l=1,2$. Therefore,\vspace*{-2pt} we have
\begin{eqnarray*}
\bigl|g_\theta(y_1) - g_\theta(y_2) \bigr| &
\leq  & d \sum_{j=1}^{m_n} h_n^\beta
\bigl|\vartheta \bigl(\rho(z_{j,n},y_1)/h_n \bigr) -
\vartheta \bigl(\rho (z_{j,n},y_2)/h_n \bigr) \bigr|
\\[-2pt]
& \leq & 2d \rho(y_1,y_2)^\beta\max \bigl\{\|
\vartheta\|_\infty,\bigl\|\vartheta'\bigr\|_\infty \bigr\},
\end{eqnarray*}
so that a sufficiently small choice of $d$ guarantees that
$g_\theta\in{\mathcal G}$ uniformly for all $\theta\in\{0,1\}^{m_n}$.

Now we use Assouad's lemma, which is based on the common Bayesian
approach of imposing the uniform distribution on $\{0,1\}^{m_n}$ to be
the a-priori distribution of $\theta$. We refer to the book of Tsybakov \cite{r23}, in particular, Section~2.7.2 for a detailed review and proof
of these results. From there, it follows\vspace*{-2pt} that
\begin{eqnarray*}
&& \hspace*{-4pt}\sup_{g \in{\mathcal G}} E_g \bigl\|\hat{g} - g\bigr\|_{P_X}^2
\\[-2pt]
&&\hspace*{-4pt}\quad\geq  \frac{1}4 d^2 h_n^{2\beta} \sum
_{j=1}^{m_n} \int_{B_{j,n}}
\vartheta^2 \bigl(\rho(z_{j,n},x)/h_n \bigr)
\,\mathrm{d}P_X(x)
\bigl\{1 - E H^2 \bigl(E_\theta f_{\theta,j,0}(y
\vert{\mathbf X}_n),E_\theta f_{\theta
,j,1}(y\vert{\mathbf
X}_n) \bigr) \bigr\},
\end{eqnarray*}
where $H^2(f_1,f_2) := \int(\sqrt{f_1}-\sqrt{f_2})^2$
denotes the squared Hellinger distance between two densities $f_1$ and
$f_2$. We consider that
\begin{eqnarray*}
&& \int_{B_{j,n}} \vartheta^2 \bigl(
\rho(z_{j,n},x)/h_n \bigr) \,\mathrm{d}P_X(x)  =
\vartheta^2(0) m_n^{-1}.
\end{eqnarray*}
Therefore, we realize that the uniform squared risk is
bounded from below by a global constant times $h_n^{2\beta}$ whenever
we can show that
\begin{equation}
\label{eqcond}
\lim_{n\to\infty} \max_{j=1,\ldots,m_n} E
H^2 \bigl(E_\theta f_{\theta
,j,0}(y\vert{\mathbf
X}_n),E_\theta f_{\theta,j,1}(y\vert{\mathbf
X}_n) \bigr) = 0.
\end{equation}
By the Cauchy--Schwarz inequality with respect to $E_\theta$
we deduce that
\[
E H^2 \bigl(E_\theta f_{\theta,j,0}(y\vert{\mathbf
X}_n),E_\theta f_{\theta
,j,1}(y\vert{\mathbf
X}_n) \bigr) \leq E_\theta E H^2
\bigl(f_{\theta
,j,0}(y\vert{\mathbf X}_n),f_{\theta,j,1}(y\vert{
\mathbf X}_n) \bigr).
\]
We consider that
\begin{eqnarray*}
&& H^2 \bigl(f_{\theta,j,0}(y\vert{\mathbf X}_n),f_{\theta
,j,1}(y
\vert{\mathbf X}_n) \bigr)
\\[-2pt]
&&\quad = 2 - 2 \prod_{k=1}^n \biggl(1 -
\frac{1}2 H^2 \bigl(f_\varepsilon \bigl(
\cdot-g_{\theta,j,0}(X_k) \bigr),f_\varepsilon \bigl(
\cdot-g_{\theta,j,1}(X_k) \bigr) \bigr) \biggr),
\end{eqnarray*}
almost surely. Applying the expectation yields that
\begin{eqnarray*}
&& E  H^2 \bigl(E_\theta f_{\theta,j,0}(y\vert{\mathbf
X}_n),E_\theta f_{\theta,j,1}(y\vert{\mathbf
X}_n) \bigr)
\\
&&\quad \leq 2 - 2 E_\theta \biggl(1 - \frac{1}2 E H^2
\bigl(f_\varepsilon \bigl(\cdot-g_{\theta
,j,0}(X_1)
\bigr),f_\varepsilon \bigl(\cdot-g_{\theta,j,1}(X_1) \bigr)
\bigr) \biggr)^n
\\
&&\quad \leq 2 - 2 E_\theta \biggl(1 - \frac{1}8 E\bigl|g_{\theta,j,1}(X_1)
- g_{\theta,j,0}(X_1)\bigr|^2\cdot\int\bigl|f_\varepsilon'(t)\bigr|^2
/ f_\varepsilon (t)\, \mathrm{d}t \biggr)^n,
\end{eqnarray*}
where, for all $\theta\in\{0,1\}^{m_n}$ and $j=1,\ldots
,m_n$, we have
\[
E\bigl|g_{\theta,j,1}(X_1) - g_{\theta,j,0}(X_1)\bigr|^2
= d^2 h_n^{2\beta} \int_{B_{j,n}}
\vartheta^2 \bigl(\rho(z_{j,n},x)/h_n \bigr)
\,\mathrm{d}P_X(x) = d^2 h_n^{2\beta}
\vartheta^2(0) m_n^{-1}.
\]
Therefore, recalling that $h_n = \delta_n/4$ we put $\delta_n
= \{c_h \log n\}^{-1/\gamma}$ for some $c_h> 1/c_{x,0}$ so that
$h_n^{2\beta} m_n^{-1} = \mathrm{o}(1/n)$ and (\ref{eqcond}) is fulfilled. This
provides the desired lower bound.

(b)  We take $m_n$, the $z_{1,n},\ldots,z_{m_n,n} \in{\mathcal
Y}$ and the balls $B_{j,n}$ from the proof of part (a). As the
$B_{j,n}$ are pairwise disjoint we have that
\[
\sum_{j=1}^{m_n} P_{X,n}
(B_{j,n} ) = P_{X,n} \Biggl(\bigcup
_{j=1}^{m_n} B_{j,n} \Biggr) \leq1,
\]
so that, for at least one $k_n \in\{1,\ldots,m_n\}$, we have
$P_{X,n} (B_{k_n,n} ) \leq1/m_n$ where $z_{k_n,n} \in{\mathcal Y}$.
To simplify the notation, we write $w_n := z_{k_n,n}$.

We consider the mappings $g_0 \dvt\equiv0$ and $g_n(z) := d h_n^\beta
\vartheta(\rho(w_n,z) / h_n)$ with $h_n = \delta_n/4$ on the domain
${\mathcal X}$ with $\vartheta$ as in the proof of (a). Again, choosing $d$
small enough ensures that $g_0, g_n \in{\mathcal G}$ for all $n$.

We define
\[
\alpha_n := \bigl|g_n(w_n) - g_0(w_n)\bigr|/2
= d \vartheta(0) h_n^\beta/2,
\]
and the events
\[
H_n(g) := \bigl\{\omega\in({\mathcal X}\times\mathbb{R})^n
\dvt \bigl|\hat{g}_n(w_n,\omega) - g(w_n)\bigr|\geq
\alpha_n \bigr\}.
\]
Also we write ${\mathbf X}_n := (X_1,\ldots,X_n)$. We deduce that
\begin{eqnarray}
\nonumber
&& \sup_{g\in{\mathcal G}} \sup_{y \in{\mathcal Y}}
P_g \bigl[\bigl|\hat{g}_n(y,{\mathcal Z}_n) -
g(y)\bigr|> d \vartheta(0)  h_n^\beta/2 \bigr]\\
\nonumber
  &&\quad\geq   \sup
_{g\in
{\mathcal G}} E P_g \bigl(H(g) \vert{\mathbf
X}_n \bigr)
\nonumber
\\[-8pt]
\label{eqnneu1}
\\[-8pt]
\nonumber
&&\quad \geq   \frac{1}2\cdot E \bigl\{P_{g_0}
\bigl(H(g_0)\vert{\mathbf X}_n \bigr) + P_{g_n}
\bigl(H(g_n)\vert{\mathbf X}_n \bigr) \bigr\}
\\
\nonumber
 &&\quad \geq   \frac{1}2 - \frac{1}2\cdot E \operatorname{TV}
(P_{g_0}\vert{\mathbf X}_n,P_{g_n}\vert{\mathbf
X}_n ),
\end{eqnarray}
where $\operatorname{TV}(P,Q)$ denotes the total variation distance
between some probability measures $P$ and~$Q$. Note that the
conditional probability measure $P_g\vert{\mathbf X}_n$ just turns out to be
the probability measure of independent random variables $\delta_j$,
$j=1,\ldots,n$, with the density $f_\varepsilon(\cdot- g(X_j))$
conditionally on ${\mathbf X}_n$. By LeCam's inequality, we have
\[
\operatorname{TV} (P_{g_0}\vert{\mathbf X}_n,P_{g_n}\vert{
\mathbf X}_n ) \leq \Biggl\{1 - \prod_{j=1}^n
\biggl(1 - \frac{1}2 H^2 \bigl(f_\varepsilon \bigl(\cdot
-g_0(X_j) \bigr),f_\varepsilon \bigl(
\cdot-g_n(X_j) \bigr) \bigr) \biggr)^2 \Biggr
\}^{1/2},
\]
almost surely, where $H(f_1,f_2)$ denotes the Hellinger
distance between two densities $f_1$ and~$f_2$. Applying the
expectation to both sides, Jensen's inequality and some information
theoretic arguments yield that
\[
E \operatorname{TV} (P_{g_0}\vert{\mathbf X}_n,P_{g_n}\vert{
\mathbf X}_n ) \leq \biggl\{1 - \biggl(1 - \frac{1}8 E
\bigl(g_n(X_1) - g_0(X_1)
\bigr)^2 \int \bigl|f_\varepsilon'(x)\bigr|^2 /
f_\varepsilon(x) \,\mathrm{d}x \biggr)^{2n} \biggr\}^{1/2}.
\]
Since the restrictions of $g_n$ and $g_0$ to the domain
${\mathcal Y}$ coincide on ${\mathcal Y} \setminus B_{\mathcal Y}(w_n,h_n)$ we
deduce that
\[
E \bigl(g_n(X_1) - g_0(X_1)
\bigr)^2 \leq d^2 h_n^{2\beta} \|
\vartheta \|_\infty^2 P_{X,n} \bigl(B_{\mathcal Y}(w_n,h_n)
\bigr) \leq d^2 h_n^{2\beta} m_n^{-1}
\|\vartheta\|_\infty^2,
\]
so that (\ref{eqnneu1}) is bounded away from zero whenever
$h_n^{2\beta} m_n^{-1}   =   {\mathcal O}(1/n)$. Under the selection of
$h_n$ and $\delta_n$ as in part (a) this condition is satisfied. That
completes the proof.
\end{pf*}

\begin{pf*}{Proof of Theorem~\protect\ref{T4}}
The inequality (\ref{eq2.3}) yields that the set ${\mathcal Y}$ contains infinitely many
elements. Fix three different $y_1,y_2,y_3 \in{\mathcal Y}$. We introduce
the sets ${\mathcal Y}_j$, $j=1,2,3$, with
\[
{\mathcal Y}_j := \bigl\{y\in{\mathcal Y} \dvt
\rho(y,y_j) \leq\rho (y,y_k), \forall k \in\{1,2,3\}
\bigr\},
\]
whose union includes ${\mathcal Y}$ as a subset. Note that
\begin{equation}
\label{eqT.4.1}
\tfrac{1}3 {\mathcal N}_{\mathcal X}(\delta ,{\mathcal
Y},\rho) \leq{\mathcal N}_{\mathcal X}(\delta,{\mathcal Y}_j,\rho)
\leq {\mathcal N}_{\mathcal X}(\delta,{\mathcal Y},\rho),
\end{equation}
holds true for all $\delta>0$ for at least one $j=1,2,3$.
Select $j=1,2,3$ such that ${\mathcal Y}':={\mathcal Y}_j$ satisfies the above
inequality; and put $z_{-1}$, $z_0$ equal to the other $y_k$, $k\neq
j$. Note that
\[
\rho(y,z_{l}) \geq M := \min \bigl\{\rho(y_r,y_s)
\dvt r\neq s \bigr\} /2,
\]
holds for all $y\in{\mathcal Y}'$ and $l=0,-1$. Clearly, we have
$\rho(z_0,z_{-1}) \geq M$ as well. The inequalities (\ref{eq2.1}),
(\ref{eq2.2}), (\ref{eq2.3}) and (\ref{eqT.4.1}) yield the existence
of $z_{1},\ldots,z_{d_n} \in{\mathcal Y}'$ with some even number $d_n \geq
\lfloor\exp(c_{x,0} \delta_{n}^{-\gamma})/3 \rfloor- 1$ such that
$\rho(z_{j},z_{k}) > \delta_{n}$ for any sequence $(\delta_{n})_n
\downarrow0$. Therefore, the balls ${\mathcal B}_{\mathcal Y}(z_{j},\delta
_{n}/4)$, $j=1,\ldots,d_{n}$ are pairwise disjoint. By $R_{n}$ we
denote the discrete probability measure which fulfills
\begin{eqnarray*}
R_n \bigl(\{z_0\} \bigr) & =& R_n \bigl(
\{z_{-1}\} \bigr) = 2\kappa M_0^{-\beta} / C,
\\
R_n \bigl(\{z_j\} \bigr) & =& \bigl(1 - 4 \kappa
M_0^{-\beta} / C \bigr) / d_n ,\qquad  j=1,
\ldots,d_n,
\end{eqnarray*}
and $M_0 := \min\{C^{-1/\beta},M\}$. We define the functions
\[
f_\theta(y) := 1 + \theta_0 \frac{1}2 C
M_0^\beta \bigl(1_{\{z_{-1}\}
}(y) - 1_{\{z_{0}\}}(y)
\bigr) + \frac{1}2 C \delta_n^\beta\sum
_{j=1}^{d_n/2} \theta_j
\bigl(1_{\{z_{2j-1}\}}(y) - 1_{\{z_{2j}\}}(y) \bigr),
\]
where $\theta:= (\theta_0,\ldots,\theta_{d_n/2})$ denotes
some binary vector. Therein we stipulate that
\begin{equation}
\label{eqkappa}
\kappa\in \bigl(0,M_0^\beta C / 8 \bigr),
\end{equation}
which does not depend on $n$. For $n$ sufficiently large,
$f_\theta$ is bounded by $1/2$ from below and by $3/2$ from above --
uniformly with respect to the vector $\theta$. Furthermore, the
functions $f_\theta$ integrate to one with respect to the probability
measure $R_n$ so that the functions $f_\theta$ are probability
densities. The probability measure generated by $f_\theta$ is denoted
by $P_\theta$.

We write $\theta'$ for the corresponding vector $\theta' := (1-\theta
_0,\ldots,1-\theta_{d_n/2})$. Then
\[
\operatorname{TV}(P_\theta,P_{\theta'}) = \frac{1}2 \int
\bigl|f_\theta(y) - f_{\theta'}(y) \bigr|\, \mathrm{d}R_n(y) \geq
\frac{1}4 \bigl(C M_0^{\beta} R_n \bigl(\{
z_{0}\} \bigr) + C M_0^{\beta} R_n
\bigl(\{z_{-1}\} \bigr) \bigr) = \kappa,
\]
so that $(P_\theta,P_{\theta'}) \in{\mathcal P}_\kappa$.
Furthermore, we have
\[
\frac{\mathrm{d}P_\theta}{\mathrm{d}(P_\theta+P_{\theta'})}(y) = \frac{f_\theta
(y)}{f_\theta(y) + f_{\theta'}(y)},
\]
so that
\begin{eqnarray*}
&& \sup_{\theta''=\theta,\theta'} \biggl|\frac{\mathrm{d}P_{\theta''}}{\mathrm{d}(P_\theta
+P_{\theta'})}(y) - \frac{\mathrm{d}P_{\theta''}}{\mathrm{d}(P_\theta+P_{\theta'})}(x) \biggr|
\leq\max \bigl\{ \bigl|f_\theta(x) - f_\theta(y) \bigr|,
\bigl|f_{\theta
'}(x) - f_{\theta'}(y) \bigr| \bigr\}.
\end{eqnarray*}
For $n$ sufficiently large (precisely, for $\delta_n <
2^{-1/\beta} M_0$), we can verify that
\[
\bigl|f_{\theta}(x) - f_{\theta}(y) \bigr| \leq C \rho(x,y)^\beta,
\]
for all $x,y \in{\mathcal Y}_0 = \{z_{-1},z_0,\ldots,z_{d_n}\}$
and all $\theta\in\{0,1\}^{d_n/2+1}$. This provides that $(P_\theta
,P_{\theta'}) \in{\mathcal P}_{C,\beta,\kappa}$ for all $\theta\in\{0,1\}
^{d_n/2+1}$ for $n$ large enough.

The underlying statistical experiment is less informative than the
model, in which exactly $n$ i.i.d. training data are drawn for each
group, that is, we observe the samples $X_1,\ldots,X_n$ from $P_X$ and
$Y_1,\ldots,Y_n$ from $P_Y$. Therefore, as we are proving a lower
bound, we may switch to the latter statistical model. As in the proof
of Theorem~\ref{T1}(a), we apply Assouad's lemma (see, e.g., Tsybakov \cite{r23}) and the Bretagnolle--Huber inequality
(see Bretagnolle and Huber \cite{r6}), which yield that
\begin{eqnarray}
\nonumber
{\mathcal E}_{n}(\varphi_n)& \geq &
\frac{1}4 \sum_{b=0}^1 \sum
_{l=0}^{d_n/2} R_n \bigl(
\{z_{2l}\} \bigr) \cdot C \bigl(M_0^\beta1_{\{0\}}(l)
+ \delta_n^\beta1_{(0,\infty)}(l) \bigr)
\\
\nonumber
&&\hspace*{19pt}\qquad{}\times \bigl(1 - E_\theta \bigl\{1 - \exp \bigl(- n {
\mathcal K}(f_{\theta,l,1},f_{\theta,l,0}) - n {\mathcal
K}(f_{\theta',l,0},f_{\theta',l,1}) \bigr) \bigr\}^{1/2} \bigr),
\end{eqnarray}
where ${\mathcal K}$ denotes the Kullback--Leibler distance
between some densities. For $l\geq1$, we deduce that
\begin{eqnarray}
\nonumber
{\mathcal K}(f_{\theta,l,1},f_{\theta,l,0}) & =& \int \biggl\{
\log\frac
{f_{\theta,l,1}(x)}{f_{\theta,l,0}(x)} \biggr\} f_{\theta,l,1}(x)\, \mathrm{d}R_n(x)
\\
\label{eqclneu4}
& \leq &  \int\frac{f_{\theta,l,1}(x)}{f_{\theta
,l,0}(x)} \bigl|f_{\theta,l,1}(x) -
f_{\theta,l,0}(x) \bigr| \,\mathrm{d}R_n(x)
\\
& \leq &  3 C \delta_n^\beta R_n
\bigl(\{z_{2l}\} \bigr) \leq 3 C \bigl(1 - 4\kappa M_0^{-\beta}
/ C \bigr) \delta_n^\beta d_n^{-1},
\nonumber
\end{eqnarray}
almost surely. The same upper bound can be established for
${\mathcal K}(f_{\theta',l,0},f_{\theta',l,1})$ analogously. Now we specify
\[
\delta_n = (c_\delta\log n)^{-1/\gamma},
\]
for some constant $c_\delta>0$. Choosing $c_\delta$
sufficiently large, (\ref{eqclneu4}) yields that
\begin{eqnarray}
\nonumber
\liminf_{n\to\infty} \delta_n^{-\beta} {
\mathcal E}_{n}(\varphi_n) & \geq  & \liminf
_{n\to\infty} \frac{1}4 C \sum_{b=0}^1
\sum_{l=1}^{d_n/2} R_n \bigl(
\{z_{2l}\} \bigr)
\\
\nonumber
& \geq &  \frac{1}4 \bigl(1 - 4\kappa M_0^{-\beta}
/ C \bigr) \geq \frac{1}8,
\end{eqnarray}
when using (\ref{eqkappa}) in the last step. The selection
of $\delta_n$ completes the proof.
\end{pf*}

\section*{Acknowledgements}
The author is grateful to the associate editor and three referees for
their inspiring comments.


\printhistory
\end{document}